\documentclass[12pt]{article}
\usepackage{epsfig}
\usepackage{color}
\usepackage{xcolor}
\usepackage{graphicx,enumerate}
\usepackage{transparent,color}
\usepackage{tikz}
\usepackage[latin1]{inputenc}
\usepackage{verbatim}
\usepackage{multicol}
\usepackage{alltt}
\usepackage{fancyvrb}
\usepackage{listings}
\usepackage{amsfonts,amsmath,amsthm,amssymb}
\usepackage{authblk}

\makeatletter
\newcommand{\leqnomode}{\tagsleft@true}
\newcommand{\reqnomode}{\tagsleft@false}
\makeatother
\begin{document}

\title{To replace or not to replace in finite population sampling}
\author{Daniel Q. Naiman \& Fred Torcaso\\ \small Department of Applied Mathematics and
Statistics\\Johns Hopkins University\\Baltimore, MD}

\maketitle

\begin{abstract}
We revisit the classical result in finite population sampling which states that in {\em
equally-likely} ``simple" random sampling the sample mean
is more reliable when we do not replace after each draw.  In this paper we investigate if
and when the same is true for samples where it may no longer be true that each member
of the population has an equal chance of being selected.
For a certain class of sampling schemes, we are able to obtain convenient expressions
for
the variance of the sample mean and surprisingly,
we find that for some selection distributions a more reliable estimate of the population
mean will happen by replacing after each draw.  We show for selection distributions lying
in a certain polytope the classical result prevails.
\end{abstract}

\noindent 
{\bf 1. Introduction.}

For samples from a finite population it is a well known, elementary, and intuitive fact
that the variance of the sample mean for a simple random sample \emph{without} replacement
does not exceed the variance of the sample mean for a sample \emph{with} replacement. See  \cite{rice},
chapter 7 for example.

Both \cite{Basu1958} and \cite{RajKhamis1958} argue that the above comparison is not necessarily the
appropriate one since sampling with replacement allows for the possibility of sampling the same element
multiple times. They go on to point out that in equal probability sampling \emph{with} replacement
the variance of the sample mean becomes more efficient if replicates are ignored.  A short proof of this
fact is given in \cite{KorwarSerfling1970}, who also provide an
expression for the reduction in variance that is amenable to numerical approximation.
A generalization of the  above results to general convex loss functions is give by \cite{SinhaSen1989}.

Unequal probability/biased sampling-based designs have a long history in sample surveys
\cite{SmithSampling}.
When probability weights are known, using the inverse-probability weighted
sample mean \cite{HT52} leads to an unbiased estimator of the population mean.
The situation gives rise to a natural question, namely, for weighted  sampling without
replacement is it still the case that the Horvitz-Thompson estimator has a
variance that is no greater than the variance of corresponding estimator
based on samples with replacement?
To address this question we need to be precise about what we mean by weighted sampling
without replacement, that is, we need to specify the joint sampling distribution. Many finite
population sampling schemes have been devised that achieve weighted sampling. In
particular, \cite{brewer2013sampling} provide a comprehensive treatment of
this subject.

We require that our sampling distribution/scheme satisfy certain natural requirements. First, the
one-dimensional marginals for the individuals drawn should be identical, as specified by
the probability weights. In addition, since the sample variance is determined from the
bivariate marginals we require that all pairwise bivariate marginals be identical as well.
Beyond these requirements, there are still potentially many choices of sampling schemes,
and investigating the problem at hand remains quite unwieldy without additional
restrictions. We choose to study a particularly simple sampling scheme wherein each bivariate
marginal probability is a particular affine function of the sum of two marginal
probabilities.

In the process of studying sampling schemes satisfying these requirements,  we found that
we could explicitly prescribe the full joint distribution with marginals of every dimension
having an affine form in the marginal probabilities, resulting in exchangeability.
There has been  considerable focus in the  literature on finite population sampling,
including
references in which sampling schemes are developed for producing
sequences with prescribed individual and pairwise inclusion probabilities \cite{bksinha73,
mukhop72}, as well as interest in finite exchangeable sequences \cite{Schervish,
diaconisfreedman80}.
Despite the simplicity of the sampling schemes we focus on, we have not found any
references to our construction. Still, we would not be the least surprised to find that it
has been introduced previously.

This paper is organized as follows.
In Section 2, we introduce some notation and formalize the problem at hand.
We focus on sampling without replacement schemes in which
the univariate marginals are identical, the bivariate marginals are identical, and we seek results that
hold for a given sampling and no matter what numerical values
associated with population individuals happen to be.
Consequently, we can associate with a sampling scheme a certain quadratic form, and formulate our
question in terms nonnegative definiteness of this form.
In Section 3 we introduce the assumption of a specific \emph{affine}
form for bivariate marginal distributions, we reexpress the quadratic form for this special case.
In Section 4, we give conditions guaranteeing that the quadratic form in question
is positive semidefinite, so that the Horvitz-Thompson estimator is guaranteed to perform at least as well
for sampling without replacement
as for sampling with replacement.
In Section 5, we give a class of sampling schemes, which we refer to as \emph{affine sampling without
replacement schemes}
in which the univariate marginals can be specified, and the bivariate marginal distributions
have the affine form required in the previous sections. We provide conditions stating that if the
univariate marginal distribution
lies in a certain polytope, then such a sampling scheme exists. 
We investigate the positive definiteness of the above-mentioned quadratic form
at the vertices of this polytope.
In Section 6, we describe algorithmically how to sample from an element of this family of distribution,
under the condition that the number of distinct marginal probabilities is small.
Section 7 provides some concluding remarks.

\noindent 
{\bf 2. Finite population sampling preliminaries}

We need to introduce a bit of notation. Assume a fixed population of size $N$ with
individuals labeled $1,2,\ldots,N.$  Each individual $i$ has an associated probability
$p_i,$ so that $p_i\geq 0$ and $\sum_{i=1}^N p_i =1.$ We use ${\bf p}$ to denote probability
vector $(p_1,\ldots,p_N).$
In addition, we assume there is a vector ${\bf x} = (x_1,\ldots,x_N)$ representing the values of
some numerical attribute for each population member.

We will use $I_1,\ldots,I_n$ to denote a sample of size $n$ drawn from the population,
where each $I_i$ has $p$ as its marginal distribution, i.e.
\begin{equation}
\label{marginalcondition}
P[ I_i = u] = p_u, \textrm{ for } u=1,\ldots,N, \textrm{ and } i=1,\ldots,n.
\end{equation}
This gives rise to a sequence of random variables $X_1,\ldots,X_n$ defined by $X_i =
x_{I_i},$ for $i=1,\ldots,n.$

Since the sampling procedure is potentially biased, we correct each $X_j$ to give the (unbiased)
Horvitz-Thompson estimator \cite{HT52}
\begin{equation}
\label{HVestimator}
\hat{\mu} = \frac{1}{n}\sum_{i=1}^n X_i/(Np_{I_i}),
\end{equation}
of the population mean $\mu = \frac{1}{N}\sum_{i=1}^N x_i.$

For all of the sampling schemes we consider, it is assumed that
the bivariate marginal distributions of the $I_i$ are identical,
that is, $\delta_{uv} = P[I_i = u,I_j = v]$ does not depend on the choice of $i\neq j.$
Then, we can express the variance of the estimator in (\ref{HVestimator}) as
\begin{equation}
\label{varianceofHV}
\textrm{Var}(\hat{\mu}) =
\frac{1}{n}\textrm{Var}(X_i/(Np_{I_i})) +
\frac{(n-1)}{n}\textrm{Cov}(X_i/(Np_{I_i}),X_j/(Np_{I_j})).
\end{equation}
Under sampling \emph{with replacement}, the covariance term in (\ref{varianceofHV}) vanishes, so we can
express the relationship between the variances under \emph{sampling without replacement} and
\emph{sampling with replacement} by
$$
\textrm{Var}(\hat{\mu}_{w/o~rep})
= \textrm{Var}(\hat{\mu}_{with~rep}) -
\frac{n-1}{n} \frac{1}{N^2}
\sum_{u,v=1}^N
(1 - \frac{\delta_{uv}}{p_up_v})x_ux_v.
$$
Consequently, under a sampling distribution with bivariate marginals given by the $\delta_{uv}$ we see
that the question of whether sampling with replacement is not better than sampling without replacement
for all possible
choices of the vector ${\bf x}$ comes down to the positive semi-definiteness of the $N \times N$
matrix
\begin{equation}
\label{criticalmatrix}
\Psi = \left( 1 -
\frac{\delta_{uv}}{p_up_v}
\right)_{1 \leq u,v \leq N}.
\end{equation}

As a special case, under the familiar uniform sampling without replacement scheme, where we have
$p_u = 1/N,$ for $u=1,\ldots,N,$ and
$$
\delta_{uv} = \left\{
\begin{array}{ll}
0 &      \mbox{ if } 1 \leq u = v \leq N, \\
-1/(N-1) & \mbox{ if } 1 \leq u \neq v \leq N,
\end{array}
\right.
$$
so that
$$
\Psi = \frac{N}{N-1}I_N + \frac{-1}{N-1}J_N
$$
where $J_N$ denotes the $N \times N$ matrix all of whose entries are 1.
For $N \geq 2,$ this matrix has eigenvalues of 0 with multiplicity 1 and $\frac{N}{N-1}$ with multiplicity $N-1$ hence is positive semi-definite.
So we conclude that for this case, sampling with replacement is never beats sampling without replacement, for
all values of ${\bf x}.$

\noindent
{\bf 3. A Class of Sampling Without Replacement Schemes}

In addition to assuming (\ref{marginalcondition}) holds,
we focus on sampling schemes whose bivariate marginals have a specific affine form, namely, where
\begin{equation}
\label{bivariatemarginalform}
\begin{array}{l}
\delta_{uv} = P[I_i=u,I_j=v] = \left\{
\begin{array}{ll}
A + B(p_u+p_v) & \textrm{ for } 1 \leq u \neq v \leq N, \\
0 &  \textrm{ for } 1 \leq u=v \leq N,\\
\end{array}
\right.\\
\textrm{ for } 1 \leq i \neq j \leq n,\\
\end{array}
\end{equation}
for some choice of constants $A$ and $B.$

The following Proposition shows that, for our purposes, the constants $A$ and $B$ are determined
in terms of $N.$

\noindent {\bf Proposition 1.}
Let $I_1,\ldots,I_n$ be random variables taking values in
$\{ 1,\ldots,N\}.$
If (\ref{marginalcondition}) and (\ref{bivariatemarginalform}) hold for constants $A$ and $B,$
and not all of the $p_u$ are equal, then $A = -1/((N-1)(N-2))$ and $B = 1/(N-2).$
Conversely, if, for some probability vector ${\bf p} = (p_1,\ldots,p_N)$ the
condition (\ref{bivariatemarginalform}) holds, where $A = -1/((N-1)(N-2)),$ and $B = 1/(N-2),$
then (\ref{marginalcondition}) holds.

{\bf Proof.}
Assume $p_u \neq p_v.$ Summing $P[I_i=u,I_j=v]$ over $u$ gives
$$
p_v = (N-1)A + B(1-p_v) + (N-1)Bp_v,
$$
and summing over $v$ gives
$$
p_u = (N-1)A + (N-1)Bp_u + (1-p_u)B.
$$
Subtracting gives
$$
p_v - p_u = -B(p_v - p_u) + (N-1)B(p_v - p_u),
$$
so that $B=1/(N-2)$ and it follows that $A = -1/((N-1)(N-2)).$

In the other direction, if (\ref{bivariatemarginalform}) holds then we have
\begin{eqnarray*}
P[I_i=u] &=& \sum_{v\neq u} P[ I_i = u, I_j = v] \\
&=& \sum_{v\neq u} A_{N,2} + B_{N,2} (p_u + p_v) \\
&=&  A_{N,2}(N-1) + B_{N,2} (N-1) p_u + B_{N,2} \sum_{v\neq u}p_v\\
&=&  A_{N,2}(N-1) + B_{N,2} (N-1) p_u + B_{N,2} (1-p_u)\\
&=&  (A_{N,2}(N-1) + B_{N,2})+ p_u((N-2)B_{N,2} - 1)= p_u.
\end{eqnarray*}
$~~\Box$

{\bf Terminology.}
For a given probability vector ${\bf p} = (p_1,\ldots,p_N)$ we will use
${\cal D}_{N,n}({\bf p})$ to denote the set of distributions
satisfying the univariate marginal condition (\ref{marginalcondition}) and
the bivariate marginal condition (\ref{bivariatemarginalform}) with
the constants $A$ and $B$ given in Proposition 1.
We may refer to these as
\emph{sampling without replacement} schemes associated with ${\bf p}$.
In contrast, we define
${\cal IID}_{N,n}({\bf p})$ to be the distribution for independent $I_1,\ldots,I_n$
satisfying (\ref{marginalcondition}) and we may refer this as the \emph{sampling with
replacement}
distribution associated with ${\bf p}.$

We address the issue of existence of sampling schemes in ${\cal D}_{N,n}({\bf p})$ in Section 5, but for now,
we note that if $p_{(1)},$ $p_{(2)}, \ldots,p_{(n)}$ denote the entries of ${\bf p}$ taken
in nondecreasing order, then by Proposition 1, the condition
\begin{equation}
\label{sumoftwosmallestprobs}
p_{(1)}+p_{(2)}\geq \frac{1}{2},
\end{equation}
guarantees the existence of a sampling distribution in ${\cal D}_{N,2}({\bf p}).$
In Section 5, we will show that for $n\geq 2$ the condition
\begin{equation}
\label{sumofnsmallestprobs}
\sum_{i=1}^n p_{(i)}\geq \frac{n-1}{N-1},
\end{equation}
is sufficient for ${\cal D}_{N,n}({\bf p})$ to be nonempty.

Having specified a probability vector ${\bf p}=(p_1,\ldots,p_N),$
we proceed to compare the variance (\ref{varianceofHV}) when $I_1,\ldots,I_n$ has a
distribution in ${\cal D}_{N,n}({\bf p})$ to the variance obtained when the distribution is
${\cal IID}_{N,n}({\bf p}).$
As indicated above, we can focus on consideration of the nonnegative
definiteness of the $N \times N$ matrix $\Psi_N(p).$
For the class of sampling schemes under
consideration,  the entries of this matrix take the form
\begin{equation}
\label{PSImatrixentries}
\psi_{ij} = \left\{ \begin{array}{ll}
1 & \mbox{ if } i = j \\
1+\left\{ \frac{1}{(N-1)(N-2)} - \frac{p_i+p_j}{(N-2)}\right\}/(p_ip_j) & \mbox{ if } i\neq
j. \\
\end{array}\right.
\end{equation}

\noindent
{\bf 4. Main results}

We proceed to sufficient conditions for nonnegative definiteness of the matrix in
(\ref{PSImatrixentries}).
To simplify matters we will make use of the following.

{\bf Lemma 1.} The eigenvalues of the matrix $\Psi = \Psi_N({\bf p})$
are all nonnegative if and only if the $(N-1)\times (N-1)$ matrix $\Gamma =\Gamma({\bf p}) =
(\gamma_{ij})$
is positive semidefinite, where
\begin{equation}
\label{gammadef}
\gamma_{ij} =
\left\{ \begin{array}{ll}
(p_i-p_N)^2 - \frac{2}{(N-1)(N-2)} + \frac{2(p_i+p_N)}{N-2} & \mbox{ for } 1 \leq i =j\leq
N-1 \\
(p_i-p_N)(p_j-p_N) - \frac{1}{(N-1)(N-2)} + \frac{2p_N}{N-2} & \mbox{ for } 1 \leq i\neq j
\leq N-1.\\
\end{array}
\right.
\end{equation}

{\bf Proof.}
Defining $y_i = x_i/p_i$ we can write ${\bf x} \Psi {\bf x}^t = {\bf y} \Omega {\bf y}^t,$
where $\Omega=\Omega({\bf p})$ denotes the $N \times N$ matrix whose $i,j$ entry is
\begin{equation}
\label{omegadef}
\omega_{ij} = \left\{ \begin{array}{ll}
p_i^2 & \mbox{ if } i = j \\
p_ip_j+\left\{ \frac{1}{(N-1)(N-2)} - \frac{p_i+p_j}{(N-2)}\right\}& \mbox{ if } i\neq j
\\
\end{array}\right.
\end{equation}
and $\Psi$ is positive semidefinite if and only if $\Omega$ is positive semidefinite.

Observe that if we take ${\bf x}={\bf p}$ the random variable $x_I/p_I$ is identically 1, so the
covariance in (\ref{varianceofHV}) vanishes.
It follows that the vector $(x_1/p_1,\ldots,x_N/p_N)$ is an eigenvector of $\Psi$ with
eigenvalue 0, and consequently
the vector ${\bf y} = (1,1,\ldots,1)$ is an eigenvector of $\Omega$ with eigenvalue 0.
Thus, to show that ${\bf y} \Omega {\bf y}^t \geq 0$ for all $y,$ it suffices to show that
$$
\sum_{i,j=1}^N \omega_{ij} y_i y_j  \geq 0
$$
whenever $\sum_{i=1}^N y_i=0,$ or equivalently, that
\begin{equation}
\label{yformnew}
\sum_{i,j=1}^{N-1} \omega_{ij} y_i y_j +
\sum_{i=1}^{N-1} \omega_{iN} y_i\left(
-\sum_{j=1}^{N-1}y_j\right) +
\sum_{j=1}^{N-1} \omega_{Nj} \left(-\sum_{i=1}^{N-1}y_i\right)y_j
+  \omega_{N,N}\left(-\sum_{j=1}^{N-1}y_j\right)^2 \geq 0,
\end{equation}
for all choices of $y_1,\ldots,y_{N-1}.$

We can rewrite expression (\ref{yformnew}) as
\begin{equation}
\label{yformnew2}
\sum_{i,j=1}^{N-1} \gamma_{ij} y_i y_j
\end{equation}
where
$$
\gamma_{ij} = \omega_{ij} - \omega_{iN} - \omega_{jN} + \omega_{NN},
$$
and it is straightforward to check that $\gamma_{ij}$ is given by (\ref{gammadef}).
$~~\Box.$

{\bf Theorem 1.} If ${\bf p}=(p_1,p_2,p_3)$ is a probability vector with
$$
p_{(1)} + p_{(2)} \geq 1/2,
$$
and
$$
\delta_{uv} = A_{N,2} + B_{N,2} (p_u+p_v), \textrm{ for } 1\leq u \neq v \leq 3,
$$
then the matrix $\Psi_N({\bf p})$ is positive semidefinite.

{\bf Proof.}
Note that $p_i \leq 1/2$ for $i=1,2.$
Using Lemma 1,  we need only show that the $2 \times 2$ matrix
$$
\Gamma = \left[
\begin{array}{cc}
(p_1-p_3)^2 - 1 + 2(p_1+p_3) & (p_1-p_3)(p_2-p_3) - 1/2 + 2p_3\\
(p_1-p_3)(p_2-p_3) - 1/2 + 2p_3 & (p_2-p_3)^2 - 1 + 2(p_2+p_3) \end{array}
\right]
$$
is positive semidefinite, and for this it suffices to show that its trace and
determinant are nonnegative.
Since $p_i+p_3\geq 1/2$ the trace is nonnegative.
Substituting $p_3 = 1- p_1-p_2,$ we obtain
$$
\det \Gamma = 18(p_1-1/2)(p_2-1/2)(p_1+p_2-1/2) \geq 0.~~~\Box
$$

{\bf Corollary 1.} If $I_1,I_2$ are random variables taking values in
$\{ 1,2,3\}$ satisfying (\ref{marginalcondition}) and (\ref{bivariatemarginalform})
where $A = -1/((N-1)(N-2)),$ and and $B=1/(N-2)$ then no matter what the value of the ${\bf x} = (x_1,x_2,x_3)$
is, the variance
of the Horvitz-Thompson estimator (\ref{HVestimator}) under a sampling without replacement
distribution in ${\cal D}_{3,2}({\bf p})$
is no greater than the variance under sampling with replacement distribution
${\cal IID}_{3,2}({\bf p}).$

Next, we turn to the case when $N>3.$

{\bf Theorem 2.} If $N>3$ and ${\bf p}=(p_1,\ldots,p_N)$ is a probability vector satisfying
$$
p_{(1)}+p_{(2)}\geq \frac{3N-2}{2N(N-1)},
$$
then the matrix $\Psi_N({\bf p})$ is positive semidefinite

{\bf Proof.}
Without loss of generality, we can assume that $p_1\geq p_2 \geq \cdots \geq p_N.$
Using Lemma 1, we can write $\Gamma({\bf p}) = R + C + D$ where  $R = R({\bf p})$ is a rank one
matrix, $C=C({\bf p})$ is a constant matrix,
and $D=D({\bf p})$ is a diagonal matrix.
$$
R = \Big( (p_i-p_N)(p_j-p_N) \Big),
$$
$$
C = \Big( -\frac{1}{(N-1)(N-2)} + \frac{2p_N}{N-2} \Big)
$$
and
$$
D = \textrm{diag}\left(-\frac{1}{(N-1)(N-2)}+\frac{2p_1}{N-2},\ldots,
-\frac{1}{(N-1)(N-2)}+\frac{2p_{N-1}}{N-2}
\right).
$$

Now let ${\bf p}^{(0)} = (1/N,\ldots,1/N).$ We can write
$$
\Gamma({\bf p}) = \Gamma({\bf p}^{(0)}) + (\Gamma({\bf p}) -\Gamma({\bf p}^{(0)}))
$$
In this case, $R({\bf p}^{(0)}) = 0,$
$C({\bf p}^{(0)})$ is the matrix whose entries are all $\frac{1}{N(N-1)}$
and $D({\bf p}^{(0)})= \frac{1}{N(N-1)}I_{N-1}.$
Thus,
$$
\Gamma(p^{(0)}) = \frac{1}{N(N-1)} (I_{N-1} + J_{N-1}).
$$
Consequently, $\Gamma({\bf p}^{(0)})$ has eigenvalues $\frac{N+1}{N(N-1)}$ with multiplicity 1,
and $\frac{1}{N(N-1)}$ with multiplicity $N-2.$

On the other hand, we have
$$
(\Gamma({\bf p}) -\Gamma({\bf p}^{(0)})) = R({\bf p}) + (C({\bf p})- C({\bf p}^{(0)}) +  
(D({\bf p})- D({\bf p}^{(0)}).
$$

We proceed to lower bound the minimal eigenvalue of each term.
\begin{itemize}
\item $R({\bf p})$ is a rank 1 matrix whose eigenvalues are given by $\sum_{i=1}^N (p_i-p_N)^2$
    with multiplicity 1, and 0 with multiplicity $N-2.$
\item $C({\bf p})-C({\bf p}^{(0)})$ is the constant matrix whose entries are all
    $\frac{2(Np_N-1)}{N(N-2)},$ so that the eigenvalues of this matrix are
$\frac{2(Np_N-1)}{N(N-2)}$  with multiplicity 1, and 0 with multiplicity $N-2.$
Since $p_N$ is the smallest of the entries of $p,$ the eigenvalue
$\frac{2(Np_N-1)}{N(N-2)}$ is non-positive, hence is minimal.
\item $D({\bf p}) - D({\bf p}^{(0)})$ is a diagonal matrix whose $i$-th diagonal entry is
    $\frac{2(Np_i-1)}{N(N-2)},$ so the eigenvalues of this matrix
are these diagonal entries and the minimal eigenvalue is $\frac{2(Np_{N-1}-1)}{N(N-2)}$
(since $p_1 \geq p_2 \geq \ldots \geq p_{N-1}.$)
\end{itemize}

Using the Rayleigh quotient for the minimal eigenvalue of a Hermitian operator \cite{Kreyszig}, 
the minimal eigenvalue of $(\Gamma({\bf p}) -\Gamma({\bf p}^{(0)}))$ is at least
$$
\frac{2(Np_N-1)}{N-2}+\frac{2(Np_{N-1}-1)}{N-2}
$$
and the minimal eigenvalue of $\Gamma({\bf p})$ is then at least
$$
\frac{1}{N(N-1)} + \frac{2(Np_N-1)}{N(N-2)}+\frac{2(Np_{N-1}-1)}{N(N-2)}
$$
and a sufficient condition
for the smallest eigenvalue to be positive is that
$$
\frac{1}{N(N-1)} + \frac{2(Np_{N}-1)}{N(N-2)}+\frac{2(Np_{N-1}-1)}{N(N-2)}>0.
$$
This is equivalent to
$$
p_{N-1} + p_{N} > \frac{3N-2}{2N(N-1)},
$$
i.e. the condition is that the sum of the two smallest $p_i$ exceeds
$\frac{3N-2}{2N(N-1)}.~~\Box$

{\bf Corollary 2.} If ${\bf p} = (p_1,\ldots,p_N)$ is a probability vector satisfying
$$
p_{(1)} + p_{(2)} \geq \frac{3N-2}{2N(N-1)},
$$
then for random variables $I_1,\ldots,I_n,$ no matter what the value of the ${\bf x} =
(x_1,\ldots,x_N)$ is, the variance of the
Horvitz-Thompson estimator (\ref{HVestimator}) under a sampling without replacement
distribution in ${\cal D}_{N,n}({\bf p})$
is no greater than the variance under sampling with replacement distribution
${\cal IID}_{N,n}({\bf p}).$

\noindent
{\bf 5. Affine Sampling Without Replacement Schemes}

The results of the previous section apply for comparing sampling without replacement schemes in ${\cal
D}_{N,n}({\bf p})$
to those in ${\cal IID}_{N,n}({\bf p})$ but so far we have not provided results concerning the
existence of elements of ${\cal D}_{N,n}({\bf p}).$ In this section, we remedy this situation and introduce a
natural generalization
of condition (\ref{bivariatemarginalform}), and Proposition 1 to higher-order marginals.

For fixed $2\leq n <N$ and a given probability vector ${\bf p}=(p_1,\ldots,p_N),$ we let 
${\cal A}_{N,n}({\bf p})$
denote the set
of joint probability distributions for random variables $I_1,\ldots,I_n$ taking values in
$\{1,\ldots,N\},$ and that
satisfy (\ref{marginalcondition}) together with
\begin{equation}
\label{affinejoint}
P[ I_1 = u_1,\ldots,I_n = u_n] =
\left\{ \begin{array}{ll}
A + B\sum_{j=1}^n p_{u_j} & u_1,\ldots,u_n \textrm{ distinct} \\
0 & \textrm{otherwise,}
\end{array}\right.
\end{equation}
for some choice of constants $A$ and $B.$
We refer to such a joint distribution as an \emph{affine sampling without replacement scheme}.

For the results that follow, we introduce constants
\begin{equation}
\label{Aform}
A_{N,k} = - \frac{(k-1)(N-k-1)!}{(N-1)!}
\end{equation}
and
\begin{equation}
\label{Bform}
B_{N,k} = \frac{(N-k-1)!}{(N-2)!}
\end{equation}
for $k=2,\ldots,N-1.$

In addition, we define constants $\tilde{A}_{N,k} = k!A_{N,k}$ and
$\tilde{B}_{N,k} = k!B_{N,k}$ for $k=2,\ldots,N-1.$
For what follows, we will need the following identities, each of which is
straightforward to verify.

\begin{equation}
\label{Identity1}
\tilde{A}_{N,n}{N \choose n} + \tilde{B}_{N,n} {N-1 \choose n-1} = 1.
\end{equation}

\begin{equation}
\label{Identity2}
\tilde{A}_{N,n}{N-1 \choose n-1} + \tilde{B}_{N,n} {N-2 \choose n-2} = 0.
\end{equation}

\begin{equation}
\label{Identity3}
\tilde{B}_{N,n}{N-2 \choose n-1} = n.
\end{equation}

\begin{equation}
\label{Identity4}
B_{N,n}(N-n) = B_{N,n-1}.
\end{equation}

\begin{equation}
\label{Identity5}
A_{N,n}(N-n+1)+B_{N,n} = A_{N,n-1}
\end{equation}

{\bf Proposition 2.} For any non-constant probability vector ${\bf p} = (p_1,\ldots,p_N),$ if a joint
distribution lies in ${\cal A}_{N,n}({\bf p}),$ so
that (\ref{affinejoint}) is satisfied for some choice of constants $A$ and $B,$ then
$A = A_{N,n},$ $B = B_{N,n},$ and
\begin{equation}
\label{inequalitycondition}
\sum_{i=1}^n p_{(i)} \geq \frac{n-1}{N-1}.
\end{equation}

{\bf Proof.}
To show that the constants $A$ and $B$ have the values claimed, we proceed by induction on $n.$
By Proposition 2, the result holds for $n=2.$
Assuming the joint distribution of $I_1,\ldots,I_n$ is in ${\cal A}_{N,n}({\bf p}),$ for distinct
values of $u_1,\ldots,u_{n-1},$  we have
$$
P[ I_1 = u_1,\ldots,I_{n-1}=u_{n-1}] = \sum_{u_n \neq u_1,\ldots,u_{n-1}} P[ I_1 = u_1,\ldots,I_n=u_n]
$$
$$
 = \sum_{u_n \neq u_1,\ldots,u_{n-1}} \left( A+ B\left\{ \sum_{j=1}^{n-1} p_{u_j}+p_{u_n} \right\}
 \right)
$$
$$
= (N-(n-1)) A + B(N-(n-1))\sum_{j=1}^{n-1} p_{u_j} + B (1-\sum_{j=1}^{n-1}p_{u_j})
$$
$$
= (A(N-(n-1)) + B) + B(N-n)\sum_{j=1}^{n-1}p_{u_j},
$$
and for non-distinct values of $u_1,\ldots,u_{n-1}$ this probability is zero. We conclude that
the joint distribution of $I_1,\ldots,I_{n-1}$ lies in ${\cal A}_{N,n-1}({\bf p}),$
so by the induction hypothesis,
$
A(N-(n-1)) + B = A_{N,n-1},
$
and
$
B(N-n) = B_{N,n-1}.
$
Using (\ref{Identity4}) and (\ref{Identity5}), we conclude that
$A = A_{N,n}$ and $B=B_{N,n}.$

To see that (\ref{inequalitycondition}) holds, observe that for distinct $u_1,\ldots,u_n$ we have
$$
A_{N,n} + B_{N,n}\sum_{j=1}^n p_{u_j}  \geq 0,
$$
that is,
$$
\sum_{j=1}^n p_{u_j}  \geq -A_{N,n}/B_{N,n} = \frac{n-1}{N-1}.
$$
$~~\Box$

For non-constant probability vectors ${\bf p}$ (the only ones we are really interested in)
we need only consider affine sampling without replacement schemes in ${\cal A}_{N,n}({\bf p})$
with the particular constants in (\ref{Aform}) and (\ref{Bform}) with $k=n.$

{\bf Theorem 3.} If ${\bf p} = (p_1,\ldots,p_N)$ is a probability vector such that
(\ref{inequalitycondition}) holds,
then taking
\begin{equation}
\label{affinejoint2}
P[ I_1 = u_1,\ldots,I_n = u_n] =
\left\{ \begin{array}{ll}
A_{N,n} + B_{N,n}\sum_{j=1}^n p_{u_j} & u_1,\ldots,u_n \textrm{ distinct} \\
0 & \textrm{otherwise,}
\end{array}\right.
\end{equation}
defines a joint probability distribution in ${\cal A}_{N,n}({\bf p}).$
In addition, random variables $I_1,\ldots,I_n$ having such a joint distribution
are exchangeable and the distribution of any $k$-tuple $I_{i_1},\ldots,I_{i_k}$
lies in ${\cal A}_{N,k}(p)$ for distinct indices $i_1,\ldots,i_k \in \{ 1,\ldots,n\}.$

Taking $k=2$ in the Theorem, we obtain as an immediate consequence the following:

{\bf Corollary 3.} ${\cal A}_{N,n}({\bf p}) \subseteq {\cal D}_{N,n}({\bf p}),$ and
consequently, assuming ${\bf p}$ satisfies (\ref{inequalitycondition})
${\cal D}_{N,n}({\bf p})\neq \emptyset.$

{\bf Proof of Theorem 3.}
We use ${\cal P}_N(n)$ to denote the set of subsets of $\{1,\ldots,N\}$ of size $n.$ 
Define
$Q:  {\cal P}_N(n) \longrightarrow \mathbb{R}$ by
$$
Q(F) = \tilde{A}_{N,n} + \tilde{B}_{N,n}\sum_{i\in F}p_i.
$$

We proceed to show
\begin{itemize}
\item[(a)] $Q(F)\geq 0$ for all $F \in {\cal P}_N(n).$
\item[(b)] $\sum_{F \in {\cal P}_N(n)}Q(F)=1,$
\item[(c)] $\sum_{F \in {\cal P}_N(n):i\in F}Q(F)=np_i,$
\end{itemize}
Once we have done this, by taking
$$
P[ I_1=u_1,\ldots,I_n = u_n] = \left\{
\begin{array}{ll}
n! Q(\{ u_1,\ldots,u_n\}) & u_1,\ldots,u_n \textrm{ distinct} \\
0 & \textrm{ otherwise},
\end{array}
\right.
$$
we obtain a distribution in ${\cal D}_{N,n}({\bf p}).$

For (a), using
(\ref{inequalitycondition}) we see that
\begin{eqnarray*}
Q(F) = \tilde{A}_{N,n} + \tilde{B}_{N,n}\sum_{i \in F}p_i &=&  -\frac {n-1}{{N-1\choose n}}
+ \frac n{{N-2\choose n-1}}\sum_{i\in F}p_i \\
& \ge & -\frac {n-1}{{N-1\choose n}} + \frac n{{N-2\choose n-1}}\frac {n-1}{N-1}
=-\frac {n-1}{{N-1\choose n}} + \frac {n-1}{{N-1\choose n}}=0.
\end{eqnarray*}

For (b) we have
\begin{eqnarray*}
\sum_{F \in {\cal P}_N(n)} Q(F) &=& \sum_{F \in {\cal P}_N(n)}
\tilde{A}_{N,n} + \tilde{B}_{N,n} \sum_{i\in F} p_i \\
&=&
\tilde{A}_{N,n} {N \choose n} + \tilde{B}_{N,n}
\sum_{F \in {\cal P}_N(n)}\sum_{i\in F} p_i\\
&=&
\tilde{A}_{N,n} {N \choose n} + \tilde{B}_{N,n}
\sum_{i=1}^N p_i \vert\{ F \in {\cal P}_N(n)~:~ i \in F\}\vert\\
&=&
\tilde{A}_{N,n} {N \choose n} + \tilde{B}_{N,n}{N-1 \choose n-1}=1.
\end{eqnarray*}
by (\ref{Identity1}).

To verify (c), for fixed $i,$ let
$$
{\cal P}_N(n,i) = \left\{ F\in {\cal P}_N(n) ~:~ i\in F\right\}.
$$
then we have
\begin{eqnarray*}
\sum_{F\in {\cal P}_N(n,i)} Q(F) &=&
\tilde{A}_{N,n} {N-1\choose n-1}
+ \tilde{B}_{N,n} \sum_{F \in {\cal P}_N(n,i)} \sum_{j \in F} p_j
\\
&=& \tilde{A}_{N,n} {N-1\choose n-1}
+\tilde{B}_{N,n} \left\{\sum_{F \in {\cal P}_N(n,i)}
\left(p_i + \sum_{j\in F:j\neq i}p_j\right) \right\}
\\
&=& \tilde{A}_{N,n} {N-1\choose n-1} + \tilde{B}_{N,n}
\left\{ p_i {N-1 \choose n-1} +  \sum_{j\neq i} p_j {N-2 \choose n-2}\right\}
\\
&=& \tilde{A}_{N,n} {N-1\choose n-1} +
\tilde{B}_{N,n} \left\{ p_i \left({N-1 \choose n-1} - {N-2 \choose n-2} \right)+
\sum_{j=1}^N p_j {N-2 \choose n-2}\right\}
\\
&=& \tilde{A}_{N,n} {N-1\choose n-1} + \tilde{B}_{N,n} {N-2 \choose n-1}p_i +
\tilde{B}_{N,n}{N-2 \choose n-2}.
\end{eqnarray*}

Using (\ref{Identity2}) and (\ref{Identity3}) we conclude that
$$
\sum_{F\in {\cal P}_N(n,i)} Q(F) = np_i
$$
so (c) holds.

This proves that (\ref{affinejoint2}) defines a sampling distribution in ${\cal A}_{N,n}({\bf p}).$
Exchangeability of this distribution is immediate.
Using exchangeability, it suffices to prove by backward induction on $k$ that if the joint
distribution
of $I_1,\ldots,I_k$ is  ${\cal A}_{N,k}({\bf p})$ for some $k>2,$ then
the joint distribution of $I_1,\ldots,I_{k-1}$ lies in  ${\cal A}_{N,k-1}({\bf p})$ and

To this end, we have
\begin{eqnarray*}
P[ I_j=u_j,j=1,\ldots,k-1] &=&
\sum_{u_k\neq u_1,\ldots,u_{k-1}} ( A_{N,k} + B_{N,k}\sum_{j=1}^kp_{u_j}) \\
&=& A_{N,k} (N-k+1)+ B_{N,k}\sum_{u_k\neq u_1,\ldots,u_{k-1}} \{
p_{u_k}+\sum_{j=1}^{k-1}p_{u_j}\}\\
&=&
A_{N,k} (N-k+1)+ B_{N,k}\{ (1-\sum_{j=1}^{k-1}p_{u_j})
+(N-k+1)\sum_{j=1}^{k-1}p_{u_j}\}\\
&=&
A_{N,k} (N-k+1)+ B_{N,k}+B_{N,k}(N-k)\sum_{j=1}^{k-1}p_{u_j}\\
&=&
A_{N,k-1} + B_{N,k-1}\sum_{j=1}^{k-1}p_{u_j},
\end{eqnarray*}
by (\ref{Identity4}) and (\ref{Identity5}).

$~~~\Box$

Next, we characterize set of probability vectors satisfying (\ref{inequalitycondition}).
Let
$$
{\cal S}_N = \left\{ {\bf p} = (p_1,\ldots,p_N) ~:~ p_i \geq 0, \sum_{i=1}^N p_i = 1 \right\}
$$
denote the $N-1$ simplex and define
$$
{\cal T}_{N,n} = \left\{ {\bf p} \in {\cal S}_N ~:~
\min_{F\in {\cal P}_N(n)}\sum_{i\in F} p_i \geq \frac{n-1}{N-1}  \right\}.
$$
$$
=\left\{ {\bf p} \in {\cal S}_N ~:~  p_{(1)} + \ldots +
p_{(n)}\geq \frac{n-1}{N-1} \right\}.
$$

In characterizing the sets ${\cal T}_{N,n}$ we will make use of concepts related to
polyhedral convex sets, and related linear programming, as covered in
\cite{bazaraashetty}.
It is easy to see that ${\cal T}_{N,n}$ a convex polytope contained in the $N-1$ simplex
with the uniform probability vector $(1/N,\ldots,1/N)$ in its interior.

{\bf Proposition 3.} ${\cal T}_{N,2} \supseteq \cdots \supseteq {\cal T}_{N,N-2}
\supseteq {\cal T}_{N,N-1}.$

{\bf Proof.}
Let ${\bf p} = (p_1,\ldots,p_N) \in {\cal S}_N$ with ordered
values $p_{(1)} \leq p_{(2)} \leq \ldots \leq p_{(N)}$ and suppose $p \in
{\cal T}_{N,n+1} \backslash {\cal T}_{N,n}.$ Then we have
\begin{equation}
\label{ineq1}
\sum_{i=1}^{n+1} p_{(i)} \geq \frac{n}{N-1},
\end{equation}
and
\begin{equation}
\label{ineq2}
\sum_{i=1}^{n} p_{(i)} < \frac{n-1}{N-1}.
\end{equation}
Then (\ref{ineq1}) and (\ref{ineq2}) gives
$$
p_{(n+1)} =  \sum_{i=1}^{n+1} p_{(i)} - \sum_{i=1}^{n} p_{(i)} > \frac{n}{N-1} -
\frac{n-1}{N-1} = \frac{1}{N-1}.
$$
It follows that $p_{(i)} >\frac{1}{N-1}$ for $i=n+1,\ldots,N,$ so we have
$$
\sum_{i=1}^N p_{(i)} = \sum_{i=1}^{n+1} p_{(i)} + \sum_{i=n+2}^{N}
p_{(i)} > \frac{n}{N-1} + \frac{N-(n+1)}{N-1} = 1,
$$
which is a contradiction. $~~~\Box$

We will use properties of the polytopes ${\cal T}_{N,n}$ below.
The following result describes the vertices for these polytopes.

{\bf Theorem 4.}
For $2\leq n <N$ the polytope ${\cal T}_{N,n}$ is $N-1$ dimensional.
If $n<N-1$ it has $2N$
vertices ${\bf p}^{(i,0)},~i=1,\ldots,N$ and ${\bf p}^{(i,\frac{1}{n})},~i=1,\ldots,N$ where
$$
p^{(i,0)}_j = \left\{ \begin{array}{ll} 0 & \mbox{ if } j=i\\ \frac{1}{N-1} & \mbox{ if } j \neq
i \\ \end{array}\right.,
$$
and
$$
p^{(i,\frac{1}{n})}_j = \left\{ \begin{array}{ll} \frac{1}{n} & \mbox{ if } j=i\\ \frac{n-1}{n(N-1)} &
\mbox{ if }
j \neq i \\ \end{array}\right.,
$$
and if $n=N-1$ there are $N$ vertices, namely ${\bf p}^{(i,0)},~i=1,\ldots,N$ so
${\cal T}_{N,N-1}$ is an $N-1$-dimensional simplex.

{\bf Proof.} See Appendix.

{\bf Theorem 5.} For $2\leq n <N$ the polytope ${\cal T}_{N,n}$ has exactly ${N \choose n}$ facets,
with one corresponding to each inequality of the form
\begin{equation}
\label{faceconstraints}
\sum_{i\in F} p_i \geq \frac{n-1}{N-1}, \textrm{ for } F \in {\cal P}_{N}(n).
\end{equation}
For $n<N-1$ the facet corresponding to $F$ has, as its vertices, the points ${\bf p}^{(i,0)}, i\in F$
and ${\bf p}^{(i,1/n)}, i \notin F.$

{\bf Proof.} See Appendix.

Having identified the facets and the vertices in each facet, we are in a position to describe the
vertex adjacencies, that is, which pairs of vertices form one-dimensional faces of ${\cal T}_{N,n}.$
We assume that $n<N-1$ since the polytope ${\cal T}_{N,N-1}$ forms a simplex, which implies that all
vertices
are adjacent. In the following, for fixed $N$ and a given $F \in {\cal P}_N(n)$ we use the notation ${\cal
F}_F$ to denote the
set of vertices in the facet corresponding to $F,$ i.e.
$$
{\cal F}_F = \{ {\bf p}^{(i,0)}, i\in F\} \cup \{ {\bf p}^{(i,1/n)}, i\notin F\}.
$$

Then we have the following.

\begin{itemize}
\item For any $i\in \{1,\ldots,N\}$ no facet contains both ${\bf p}^{(i,0)},$ and ${\bf p}^{(i,1/n)},$ and the same
    must be true of the intersection
of any facets, so ${\bf p}^{(i,0)}$ and ${\bf p}^{(i,1/n)}$ are not adjacent.
\item If $n=2,$ for indices $i\neq j$ the only facet containing both ${\bf p}^{(i,0)},$ and ${\bf p}^{(j,0)}$ is
    ${\cal F}_{\{i,j\}},$ and this facet contains other vertices, so no face (interection of facets) has
    a
    vertex set consisting of only these two points. These two vertices are not adjacent.
\item If $n=N-2,$ for indices $i\neq j$ the only facet containing both ${\bf p}^{(i,1/n)},$ and ${\bf p}^{(j,1/n)}$
    is ${\cal F}_{\{i,j\}^c},$ and this facet contains other vertices, so no face has a vertex set
    consisting of only these two points, hence these two vertices are not adjacent.
\item If $2<n<N-2$ and indices $i \neq j,$ for all $k\neq i,j,$ there exists a set $F \in {\cal P}_N(n)$
    such that
$i,j\in F$ and $k\notin F$ so
$$
\bigcap_{F\in {\cal P}_N(n)~:~ i,j\in F, k\notin F}{\cal F}_F = \{ {\bf p}^{(i,0)}, {\bf p}^{(j,0)}\},
$$
and hence ${\bf p}^{(i,0)}$ and ${\bf p}^{(j,0)}$ are adjacent.
\item If $2<n<N-2$ and indices $i \neq j,$ for all $k\neq i,j,$ there exists a set $F \in {\cal P}_N(n)$
    such that
$i,j\notin F$ and $k\in F$ so
$$
\bigcap_{F\in {\cal P}_N(n)~:~ i,j\notin F, k\in F}{\cal F}_F = \{ {\bf p}^{(i,1/n)}, {\bf p}^{(j,1/n)}\},
$$
and hence ${\bf p}^{(i,1/n)}$ and ${\bf p}^{(j,1/n)}$ are adjacent.
\item If $2\leq n \leq N-2$ and indices $i \neq j,$ for all $k\neq i,j,$ there exists a set $F \in {\cal
    P}_N(n)$ such that
$i\in F$ and $j,k\notin F,$ and there exists a set $F \in {\cal P}_N(n)$ such that
$i,k\in F$ and $j\notin F,$ so
$$
\bigcap_{F\in {\cal P}_N(n)~:~ i\in F, j\notin F}{\cal F}_F = \{ {\bf p}^{(i,0)}, {\bf p}^{(j,1/n)}\},
$$
and hence ${\bf p}^{(i,0)}$ and ${\bf p}^{(j,1/n)}$ are adjacent.
\end{itemize}

It follows that ${\cal T}_{4,2}$ has the face incidence structure of a 3-dimensional cube.
For $N=5,$ each facet has the face structure of a double pyramid.

The results presented so far tell us that as long as the probability vector ${\bf p}$ is
lies close enough to the center of the simplex ${\cal S}_N$ the matrix $\Psi_N({\bf p})$
is positive definite. It is natural, then, to explore what happens for ${\bf p}$ on the boundary
of
${\cal T}_{N,2}.$

{\bf Theorem 6.} For $N>3$ the matrix $\Psi_N({\bf p})$ is positive semidefinite
for all vertices ${\bf p}$ of ${\cal T}_{N,n}.$

{\bf Proof.}
Again we make use of the Lemma.
For a vertex of the type ${\bf p}^{(i,0)}$ without loss of generality we take $i=N.$
It is straightforward to show that the
diagonal entries of the $(N-1)\times (N-1)$ matrix $\Gamma$
are given by $d=\frac{1}{(N-1)^2}$ and the off-diagonal entry is $o=-\frac{1}{(N-1)^2(N-2)}.$
$\Gamma$ has eigenvalue $d+(N-2)o = 0$ with multiplicity 1, and
$d-o = \frac{1}{(N-1)(N-2)}$ with multiplicity $N-2.$

On the other hand, for a vertex of the type ${\bf p}^{(i,1/n)},$ again, without loss of generality, we take
$i=N.$
Then it is easy to check that
$\Psi$ is the $(N-1)\times (N-1)$ matrix $\Gamma$ whose diagonal entries
are all given by $d = \frac{N^2 + n(n-2)}{(N-1)^2n^2},$ and whose off-diagonal entries are all
$o = \frac{N^2(N-2) - n(n-2)}{(N-2)(N-1)^2n^2}.$
It is straightforward to check that the eigenvalues satisfy
$$
d+(N-2)o = \frac{N^2}{(N-1)n^2},
$$
and
$$
d-o = \frac{n-2}{(N-2)(N-1)n}.
$$
So all eigenvalues of $\Gamma$ are nonnegative.
$~~~\Box$

Armed with the evidence presented above, it is tempting to conjecture that
$\Psi_N({\bf p})$ is positive semidefinite for all ${\bf p} \in {\cal T}_{N,2}.$
However, the following result shows that this is not the case, even for certain choices of
$p$
on the boundary of ${\cal T}_{N,2}.$

{\bf Theorem 7.}
For $N>3$ taking ${\bf p} = \frac{1}{2}({\bf p}^{(1,1/2)}+{\bf p}^{(N,0)})$ the matrix $\Psi_N({\bf p})$ has
a negative eigenvalue.

{\bf Proof.}
Here
$$
{\bf p} = \left(\frac{1}{4} +
\frac{1}{2(N-1)},\frac{3}{4(N-1)},\frac{3}{4(N-1)},\ldots,\frac{3}{4(N-1)},
\frac{1}{4(N-1)}\right).
$$
We use the lemma once again and a straightforward computation shows that the $(N-1)\times
(N-1)$ matrix $\Gamma = \Gamma({\bf p})$
has entries
$$
\gamma_{ij} = \left\{
\begin{array}{ll}
\frac{N^2}{16(N-1)^2} - \frac{1}{2(N-1)(N-2)} + \frac{1}{2(N-2)} & \mbox{ if } i=j=1\\
\frac{N}{8(N-1)^2} - \frac{1}{2(N-1)(N-2)} & \mbox{ if } i=1, j\neq 1 \mbox{ or } i\neq 1,
j=1\\
\frac{1}{4(N-1)^2} - \frac{1}{2(N-1)(N-2)} + \frac{1}{2(N-1)(N-2)} & \mbox{ if } i,j\neq 1,
i\neq j\\
\frac{1}{4(N-1)^2}  & \mbox{ if } i,j\neq 1, i=j\\
\end{array}
\right.
$$
Let ${\bf v}$ denote the $N-1$ vector $(x,1,\ldots,1).$ Then ${\bf v}\Gamma $ takes the form
$(f,g,\cdots,g)$ where
$$
f=\left\{ \frac{N^2}{16(N-1)^2} - \frac{1}{2(N-1)(N-2)}+ \frac{1}{2(N-2)}\right\} x +
\left\{ \frac{N}{8(N-1)^2} - \frac{1}{2(N-1)(N-2)}\right\}(N-2)
$$
$$
=  \left\{ \frac{N^2}{16(N-1)^2} + \frac{1}{2(N-1)}\right\} x + \left\{
\frac{N(N-2)}{8(N-1)^2} - \frac{1}{2(N-1)}\right\}
$$
and
$$
g = \left\{ \frac{N}{8(N-1)^2} - \frac{1}{2(N-1)(N-2)}\right\} x +
\left\{\frac{1}{4(N-1)^2} - \frac{1}{2(N-1)(N-2)}\right\} (N-2) +
$$
$$
 \frac{1}{2(N-1)(N-2)}
$$
$$
= \left\{ \frac{N}{8(N-1)^2} - \frac{1}{2(N-1)(N-2)}\right\} x +
\left\{\frac{N-2}{4(N-1)^2} - \frac{N-3}{2(N-1)(N-2)}\right\}.
$$
and ${\bf v}$ is an eigenvector of $\Gamma$ with eigenvalue $g$ provided that $x$ satisfies
$f=gx.$

The equation $f=gx$ reduces to the quadratic  $Ax^2+Bx+C=0,$ where
$$
A = \frac{N}{8(N-1)^2} - \frac{1}{2(N-1)(N-2)},
$$
$$
B = \left\{\frac{N-2}{4(N-1)^2} - \frac{N-3}{2(N-1)(N-2)}\right\} -  \left\{
\frac{N^2}{16(N-1)^2} + \frac{1}{2(N-1)}\right\}
$$
$$
= -\frac{N^3 + 10N^2 - 40N +24}{16(N-2)(N-1)^2}
$$
and
$$
C =-\left\{ \frac{N(N-2)}{8(N-1)^2} - \frac{1}{2(N-1)}\right\}.
$$

The discriminant of this equation is given by
$$
D=B^2 - 4AC = \frac{N^6+36N^5-204N^4+336N^3-96N^2-128N+64}{256(N-2)^2(N-1)^4}
$$
which is positive for all $n\geq 4.$ Substituting the root $x = (-B- \sqrt{B^2-4AC})/(2A)$
into the expression
for $g$ and simplifying we conclude that
$$
\lambda = \frac{1}{32}\left\{
1 + \frac{8}{N-2}-\frac{3}{(N-1)^2} - \frac{2}{N-1}-\sqrt{256D}
\right\}
$$
is an eigenvalue of $\Gamma.$
This is negative provided that
$$
256D > \left( 1+ \frac{8}{N-2} - \frac{3}{(N-1)^2} - \frac{2}{N-1} \right)^2.
$$
which is equivalent to
$$
N^6+36N^5-204N^4+336N^3-96N^2-128N+64> $$
$$(N-2)^2(N-1)^4\left\{
1+ \frac{8}{N-2} - \frac{3}{(N-1)^2} - \frac{2}{N-1}
\right\}^2
$$
or in other words
$$
32N^5-188N^4+356N^3-236N^2+72N-36>0.
$$
It is easy to verify that this last polynomial is positive for $N \geq 1.~~\Box$

Finally, we are able to give explicit examples in which it is better not to replace than it is to
replace.

{\bf Corollary 3.} For $N>3$ and ${\bf p} = (N+1,3,\ldots,3,1)/(4(N-1))$
and $I_1,I_2$ are random variables taking values in
$\{ 1,\ldots,N\}$ there exists ${\bf x} = (x_1,\ldots,x_N)$  such that
the variance of the Horvitz-Thompson estimator (\ref{HVestimator}) under a sampling without
replacement distribution in ${\cal D}_{N,2}({\bf p})$ exceeds the variance under sampling with
replacement distribution
${\cal IID}_{3,2}({\bf p}).$

The probability distribution in this Corollary is rather pathological in that
$$
P[I_1=N,I_2=i] = A_{N,2} + B_{N,2} (1+3)/(4(N-1)) = 0, \mbox{ for } i=2,\ldots,N-1.
$$
This phenomenon arises due to fact that ${\bf p}$ lies on the boundary of ${\cal T}_{N,2}.$
In addition, we do not give an explicit construction of the vector ${\bf x}.$
By forming a mixture of this $p$ with the uniform distribution, with sufficiently small
weight given to the uniform distribution, by continuity we still obtain a matrix
$\Psi_N({\bf p})$ having a negative eigenvalue. At the same time, we can force all pairwise joint
probabilities to be positive. This leads to constructions of simple examples in which it is
better to replace than not to replace.

{\bf Example.}
Consider the case $N=4$ where
$$
{\bf p} = \frac{99}{100}(5/12,3/12,3/12,1/12) +\frac{1}{100}(1/4,1/4,1/4,1/4) =
(.415,.25,.25,.085).
$$
Here, $p_3+p_4 = .335 > 1/3$ so we can define a bivariate sampling distribution
without replacement distribution as in Proposition 1. The probability mass function of this
distribution takes the form
$$
\frac{1}{1200}\left[
\begin{array}{cccc}
0      & 199 & 199   & 100 \\
199 &     0  & 100   & 1 \\
199 &   100  & 0     & 1\\
100 &     1  & 1     & 0 \\
\end{array}
\right].
$$
The negative eigenvalue of this matrix is $\frac{4}{3} - 22\frac{\sqrt{15890}}{1411}
\approx -.6321.$ Taking ${\bf x}$ to be the corresponding eigenvector, i.e. ${\bf x} \approx
(.441,-.536,-.536,1)$ we find that for a sample $I_1,I_2$ according to the above
distribution,
then Horvitz-Thompson estimator
$$\hat{\mu}_{w/o~rep} = \frac{1}{2}\left(X_{I_1}/(4p_{I_1})+X_{I_2}/(4p_{I_2})\right)$$
we obtain
$$
\textrm{Var}(\hat{\mu}_{w/o~rep}) = .485.
$$
On the other hand, if we sample $I_1,I_2$ with replacement according to
${\bf p},$ the Horvitz-Thompson estimator has
$$
\textrm{Var}(\hat{\mu}_{with~rep}) = .450.
$$

Using the same probability vector ${\bf p}$ and taking ${\bf x}$ to be the binary vector $(1,0,0,1)$ we
obtain
$$
\textrm{Var}(\hat{\mu}_{w/o~rep}) = .341.
$$
On the other hand, if we sample $I_1,I_2$ with replacement according to
${\bf p},$ the Horvitz-Thompson estimator has
$$
\textrm{Var}(\hat{\mu}_{with~rep}) = .318.
$$

\noindent
{\bf 6. Sampling Issues}

The results in previous sections refer to affine sampling without replacement distributions, but without consideration as to how to actually implement sampling from such a distribution. For large values of the population size $N,$ and moderate values of the sample size $n,$ this can be a challenging undertaking for a general choice of the probability vector ${\bf p} = (p_1,\ldots,p_N).$ On the other hand, under the condition that the population is stratified, with individuals in each stratum assigned the same sampling probability, so that the number of distinct probabilities $p_i$ is small, and $n\ll N,$ this becomes manageable.

Assume we have $K$ distinct values of $p_i$ denoted by $p^{(1)},\ldots,p^{(K)},$ and that $p^{(i)}$ appears $N_i$ times in ${\bf p}.$ We can further assume that the population labels have been reordered so that the vector of probabilities is given by
$$
{\bf p} = (\underbrace{p^{(1)},\cdots,p^{(1)}}_{N_1},\underbrace{p^{(2)},\cdots,p^{(2)}}_{N_2},\ldots,
\underbrace{p^{(K)},\cdots,p^{(K)}}_{N_K}).
$$
Define the $j$-th stratum
$$
{\cal G}^{(j)} = \{ 1\leq i\leq N ~:~ p_i = p^{(j)}\} = \{ N_1+\ldots+N_{j-1}+1,\ldots,
N_1+\ldots+N_j \}.
$$

Fix a sample size $n,$ and assume $n\leq \min\{N_1,\ldots,N_K\}.$
Now suppose $I_1,\ldots,I_n$ has the affine sampling without replacement distribution associated with
${\bf p},$ and let $M_j$ denote the number of $I_i$ falling in ${\cal G}^{(j)},$ for $j=1,\ldots,K.$
Then for distinct indices $i_1,\ldots,i_n$ from
$\{ 1,\ldots,N\}$ we have
$$
P\left[ I_1 = i_1, \ldots,I_n =i_n\right] = A_{N,n} + B_{N,n} \sum_{j=1}^K m_j p^{(j)},
$$
where $m_j$ denotes the number of $i_p$ falling in ${\cal G}^{(j)}.$
Furthermore, the number of choices of $i_1,\ldots,i_n$ giving rise to $m_1,\ldots,m_K$ is
$n!\prod_{j=1}^k {N_j \choose m_j}.$
Consequently, ${\bf M} = (M_1,\ldots,M_K)$ has probability mass function
$$
f({\bf m}) =
n! \prod_{j=1}^K {N_j \choose m_j} \left\{A_{N,n} + B_{N,n}\sum_{j=1}^K m_jp^{(j)}\right\}
$$
where ${\bf m} = (m_1,\ldots,m_K)$ is a nonnegative integer vector with $\sum_{j=1}^K m_j = n.$
In the following, ${\bf m}$ will always denote such a vector.

By symmetry, conditionally, given ${\bf M} = {\bf m}$ the $I_1,\ldots,I_n$
are obtained by drawing a simple random sample of size $m_j$ from each stratum ${\cal G}^{(j)},$ and then randomly permuting the $n$ values obtained. Consequently, a practical sampling scheme becomes available once we are able to sample ${\bf M}.$ For this purpose, we can use an importance sampling algorithm (\cite{Fishman}).
Let
$$
g({\bf m}) = {n \choose m_1,\ldots,m_K} \prod_{j=1}^K \lambda_j^{m_j},
$$
denote the $\textrm{Multinomial}(N,\lambda_1,\ldots,\lambda_K)$
probability mass function, where $\lambda_j = N_j/N$ for $j=1,\ldots,K.$
This distribution is straightforward to sample from when $K$ is small.

Let $h({\bf m}) = f({\bf m})/g({\bf m}).$
We proceed to find an upper bound $C$ for $h({\bf m})$ depending on the choice of $n,$ and $p^{(j)}, N_j,$
for $j=1,\ldots,K.$ Then the following algorithm yields a sample
${\bf M}$ from $f_{\bf M},$ with expected number of iterations bounded by $C$:

\begin{lstlisting}
   Repeat
     Generate ${\bf  M} \sim \textrm{Multinomial} (N, \lambda_1,\ldots,\lambda_K)$
     Generate $U \sim \textrm{Uniform(0,1)}$
   Until
     $f({\bf M}) \leq CU g({\bf M})$
   Return ${\bf M}$
\end{lstlisting}

For this upper bound, we will need the following.

{\bf Lemma 2.}
For positive integers $r$ and $s$ with $1\leq r\leq s$ we have
$$
s^{-r} e^{
-\frac{r^3}{s(s-r)} -\frac{1}{2}\frac{r}{s-r} - \frac{1}{144s^2}
}\leq
\frac{(s-r)!}{s!} \leq s^{-r} e^{\frac{r^2 + \frac{1}{12}}{s-r}-\frac{1}{2}\frac{r}{s}}.
$$

{\bf Proof.} See Appendix.

{\bf Theorem 8.} The expression
$$
C =
e^{\sum_{j=1}^K \left\{
\frac{n^3}{N_j(N_j-n)} + \frac{n}{N_j-n} + \frac{1}{144N_j^2}
\right\}}
e^{\left\{\frac{n^2 + \frac{1}{12}}{N-n-2} +\frac{3}{2}\frac{n}{N-2}\right\}}
N \left( n\max_j p^{(j)} - \frac{n-1}{N-1}\right)
$$
is an upper bound for $h({\bf m}).$

{\bf Proof.}
Since $A_{N,n} = -\frac{n-1}{N-1}B_{N,n},$ we can write
$$
h({\bf m}) = \prod_{j=1}^4 h_j({\bf m}),
$$
where $h_1({\bf m}) = \prod_{j=1}^k (\frac{N}{N_j})^{m_j},$
$h_2({\bf m}) = \prod_{j=1}^K \frac{N_j!}{(N_j - m_j)!},$
$h_3({\bf m})  \equiv \frac{(N-n-2)!}{(N-2)!},$
and
$$
h_4({\bf m}) = (N-n-1)\left[  - \frac{n-1}{N-1} + \sum_{j=1}^K m_j p^{(j)}\right].
$$

Applying the upper bound in Lemma 2, we have
$$
h_2({\bf m}) \leq \left\{ \prod_{j=1}^K N_j^{m_j} \right\}
e^{\sum_{j=1}^K \left\{
\frac{m_j^3}{N_j(N_j-m_j)} + \frac{m_j}{N_j-m_j} + \frac{1}{144N_j^2}
\right\}
}.
$$
Using the fact that $m_j \leq n$ we obtain
$$
h_2({\bf m}) \leq \left\{ \prod_{j=1}^K N_j^{m_j} \right\}
e^{\sum_{j=1}^K \left\{
\frac{n^3}{N_j(N_j-n)} + \frac{n}{N_j-n} + \frac{1}{144N_j^2}
\right\}
}.
$$

Using the lower bound in Lemma 2 with $s=N-2$ and $r = n$ gives
$$
h_3 = \frac{(N-n-2)!}{(N-2)!} \leq (N-2)^{-n} e^{\frac{n^2 + \frac{1}{12}}{N-n-2} -\frac{1}{2}\frac{n}{N-2}},
$$
and using $\log(1+x) \leq \frac{x}{1+x}$ for $x>-1$ we see that
$$
(N-2)^{-n} = N^{-n} \left(\frac{N-2}{N}\right)^{-n} = N^{-n} e^{-n \log(1-\frac{2}{N})}
\leq N^{-n} e^{
-n
\frac{\frac{-2}{N}}{1 - \frac{2}{N}}
}
= N^{-n}e^{\frac{2n}{N-2}},
$$
so that
$$
h_3 \leq N^{-n} e^{\frac{n^2 + \frac{1}{12}}{N-n-2} +\frac{3}{2}\frac{n}{N-2}}.
$$
Finally, since $\sum_{j=1}^K m_j = n$ we have
$$
\sum_{j=1}^K m_j p^{(j)} = n \sum_{j=1}^K \frac{m_j}{n} p^{(j)} \leq  n\max_{j}p^{(j)},
$$
and hence
$$
h_4({\bf m}) \leq N \left[  - \frac{n-1}{N-1} + n\max_j p^{(j)}\right].
$$
The result then follows. $~~\Box$

The bound in Theorem 8 can be used to demonstrate the practicality of the importance sampling algorithm
for a host of situations. If we assume the number of strata $K$ is relatively small, that $N_j/N \approx K$
and $n\ll N_j$ for $j=1,\ldots,K,$ then the dominant term in the exponent of the expression for $C$ is the $\frac{n^2+\frac{1}{12}}{N-n-2} \approx \frac{n^2}{N}$ term. In addition, if we let $\omega = \max_j p^{(j)}/\min_j p^{(j)},$ then since
$$N \min p^{(j)} = \sum_{i=1}^N  \min p^{(j)} \leq \sum_{i=1}^N p_i = 1,$$
we have
$$
N \left( n\max_j p^{(j)} - \frac{n-1}{N-1}\right) \leq Nn\max_j p^{(j)} = Nn\omega \min_j p^{(j)} \leq n\omega.
$$
Consequently, in this situation, an approximate upper bound for the number of iterations required by the algorithm is
$n \omega e^{\frac{n^2}{N}}.$ Thus, and we conclude that the algorithm can be quite practical as long as $\frac{n^2}{N}$ is not too large.

{\bf 7. Conclusions.}

Our aim has been to compare the performance of sampling without and sampling with replacement.
In order to make a fair comparison, we assumed that both sampling procedures have the same marginal distributions,
and to go any further, we needed to make an assumption that the bivariate joint distributions in our sampling without replacement scheme have an affine structure.
For a given sample size $n,$ under the assumption that sum of the smallest marginal probabilities is sufficiently large, we constructed a  sampling scheme for sampling without replacement satisfying these conditions, as well as exchangeability, and we demonstrated that it can practical to sample from  such a distribution under suitable stratification assumptions.

It is important to ask how one might generalize our results. Specifying an exchangeable sampling distribution for a sample of size $n$ without replacement from a population of size $N$ amounts to specifying ${N \choose n}$ probabilities, which are constrained to sum to 1. The marginal probability constraint gives $N$ linear equations that these probabilities satisfy. Thus, the set of distributions one \emph{could} investigate can be viewed as a very high dimensional manifold, while the family of distributions we focus on is but one particularly mathematically convenient one, and ours can be seen as constrained to lie in a neighborhood of the uniform sampling without replacement distribution. It should not be too surprising to show that there are neighborhoods of the uniform sampling scheme in which sampling without replacement beats sampling with replacement, since this follows by a continuity argument. However, finding such neighborhoods explicitly can be a challenging undertaking. Future efforts will focus on finding other families of sampling schemes in which this can be done.

{\bf Acknowledgements.}
The authors gratefully acknowledge their colleagues James Fill and Laurent Younes for comments leading to an improved version of this manuscript.

\newpage
\begin{center}{\bf Appendix.}\end{center}

{\bf Proof of Theorem 4.}

To determine the vertices of ${\cal T}_{N,n},$ for $n<N$ for each permutation ${\bf \pi} = (\pi_1,\ldots,\pi_N)$
of $(1,\ldots,N)$ define the polyhedron
$$
{\cal O}_{\bf \pi} = \left\{ {\bf x} = (x_1,\ldots,x_N) ~:~ x_{\pi_1} \leq x_{\pi_2} \leq \cdots \leq x_{\pi_N} \right\}.
$$
Then we write ${\cal T}_{N,n}$ as a union of $N!$ polytopes
$$
{\cal T}_{N,n} = \bigcup_{{\bf \pi} \in {\cal S}_n} {\cal O}_{\bf \pi} \cap {\cal T}_{N,n}.
$$
Every vertex of ${\cal T}_{N,n}$ is a vertex of at least one of the 
${\cal O}_{\bf \pi} \cap {\cal T}_{N,n}.$

Our strategy for determining all of the vertices of ${\cal T}_{N,n}$ is to find the
vertices of each ${\cal O}_{\bf \pi} \cap {\cal T}_{N,n},$ then determine which of these
remains a vertex of ${\cal T}_{N,n}.$
By symmetry, to determine the vertices of a particular ${\cal O}_{\bf \pi} \cap {\cal T}_{N,n},$ it suffices to
consider the case when
${\bf \pi} = (1,\ldots,N).$
So we focus attention on
$$
{\cal U}_{N,n}= {\cal O}_{(1,\ldots,N)} \cap {\cal T}_{N,n}
= \left\{ {\bf x} ~:~ 0 \leq x_1 \leq \cdots \leq x_N, \sum_{i=1}^n x_i \geq \frac{n-1}{N-1},
\sum_{i=1}^N x_i = 1\right\}.
$$

{\bf Lemma.} For $2\leq n < N$ the ${\cal U}_{N,n}$
is an $N-1$-dimensional polytope
with $N$ vertices, namely, the points
\begin{equation}
\label{vertexform}
{\bf w}^{(j)} = (\underbrace{a_j,\ldots,a_j}_{j}, \underbrace{b_j,\ldots,b_j}_{N-j}), \textrm{ for }
j=1,\ldots,N,
\end{equation}
where
$$
a_j = \left\{
\begin{array}{ll}
\frac{j-1}{j(N-1)} & \textrm{ if } 1 \leq j \leq n\\
\frac{n-1}{n(N-1)} & \textrm{ if } n\leq j \leq N-1\\
\frac{1}{N} & \textrm{ if } j=N\\
\end{array}
\right.,
$$
and
$$
b_j = \left\{
\begin{array}{ll}
\frac{1}{N-1} & \textrm{ if } 1 \leq j \leq n\\
\frac{n(N-1)-j(n-1)}{n(N-j)(N-1)}& \textrm{ if } n\leq j \leq N-1\\
\end{array}.
\right.
$$

{\bf Remark 1.} For ${\bf w}^{(j)}$ of the form in (\ref{vertexform}), for $1\leq j \leq N-1,$ the values of
$a_j$ and $b_j$
are determined from the requirement that
$$
\sum_{i=1}^n w^{(j)}_i = \frac{n-1}{N-1},
$$
and
$$
\sum_{i=1}^N w^{(j)}_i = 1.
$$

{\bf Proof.}
Substituting $x_N = 1 - \sum_{i=1}^{N-1}x_i$ we can identify ${\cal U}_{N,n}$ with a polytope in
$\mathbb{R}^{N-1}$
defined by the following $N+1$ inequality constraints
$$
\begin{array}{ll}
(I_1) & \langle {\bf e}^{(1)}, {\bf x}\rangle \geq 0\\
(I_2) & \langle {\bf e}^{(2)}- {\bf e}^{(1)},{\bf x} \rangle \geq 0\\
(I_3) & \langle {\bf e}^{(3)}- {\bf e}^{(2)},{|bf x} \rangle \geq 0\\
\vdots & \vdots \\
(I_{N-2}) & \langle {\bf e}^{(N-2)}- {\bf e}^{(N-3)},{\bf x} \rangle \geq 0\\
(I_{N-1}) & \langle {\bf e}^{(N-1)}- {\bf e}^{(N-2)},{\bf x} \rangle \geq 0\\
(I_N) &  \langle -{\bf e}^{(1)}-{\bf e}^{(2)}-\cdots-
{\bf e}^{(N-2)}-2{\bf e}^{(N-1)},{\bf x} \rangle \geq -1\\
(I_{N+1}) & \langle {\bf e}^{(1)}+{\bf e}^{(2)}+\cdots+{\bf e}^{(n)},{\bf x}\rangle \geq \frac{n-1}{N-1}.\\
\end{array}
$$

Let ${\bf v}^{(j)}$ denote the $N-1$-vector, and $c_j$ the scalar, so that constraint $I_j$ takes the form
$$
(I_j) ~~~~~\langle {\bf v}^{(j)}, {\bf x} \rangle \geq c_j,
$$
for $j=1,\ldots,N+1.$
Given a point ${\bf x}  = (x_1,\ldots,x_{N-1}) \in {\cal U}_{N,n}$ define
$$
J({\bf x}) = \left\{ 1\leq j\leq N+1~:~ \langle {\bf v}^{(j)}, {\bf x} \rangle = c_j\right\}
$$
then ${\bf x}$ is a vertex of  ${\cal U}_{N,n}$ if and only if $J({\bf x}) \neq \emptyset$ and
the matrix whose rows are ${\bf v}^{(j)}, j\in J(x)$ is of rank $N-1$ (\cite{bazaraashetty}).
As a consequence, for ${\bf x}$ to be a vertex it must be the case
that $\vert J({\bf x}) \vert \geq N-1.$

We proceed to determine necessary conditions that are satisfied for a vertex ${\bf x}$
on a case by case basis. In the following, $x_N$ refers to $1-\sum_{j=1}^{N-1} x_j.$

{\bf Case 1.} Assume $1 \in J({\bf x}),$ i.e. $x_1 = 0.$
Then we have
$$
x_2 + x_3+\cdots+x_n \geq \frac{n-1}{N-1}
$$
which together with $x_2 \leq x_3 \leq \cdots \leq x_n$ implies $x_n \geq \frac{1}{N-1}.$
On the other hand
$$
x_{n+1} + x_{n+2}+\cdots+x_{N-1}+x_N = 1 - \left( x_2 + x_3 + \cdots + x_n \right)\leq 1 - \frac{n-1}{N-1}
= \frac{N-n}{N-1},
$$
which together with $x_{n+1} \leq x_{n+2} \leq \cdots \leq x_N$ implies $x_{n+1} \leq \frac{1}{N-1}.$
Since  $\frac{1}{N-1} \leq x_n \leq x_{n+1} \leq \frac{1}{N-1}$ we conclude that $x_n = x_{n+1} =
\frac{1}{N-1},$
and, in addition,
$$
x_3+\cdots+x_n \geq \frac{n-2}{N-1}
$$
and
$$
x_{n+2}+\cdots+x_{N-1}+x_N \leq \frac{N-n-1}{N-1}.
$$
Proceeding inductively, we obtain ${\bf x} = (0,\frac{1}{N-1},\frac{1}{N-1},\cdots,\frac{1}{N-1}) = {\bf w}^{(1)}.$

{\bf Case 2.} $1 \notin J({\bf x})$ and $2 \notin J({\bf x})$ 
then since $\vert J({\bf x}) \vert \geq N-1$ we must have $j
\in J({\bf x})$ for
$j=3,4,\ldots,N+1.$ Thus, $0 < x_1 < x_2 = x_3 = x_4 = \cdots = x_{N-1} = x_N.$
Now
$$
\sum_{i=1}^N x_i = x_1 + (N-1)x_2 = 1,
$$
and
$$
\sum_{i=1}^n x_i = x_1 + (n-1)x_2 =  \frac{n-1}{N-1}.
$$
Taking the difference, we see that
$$
(N-n)x_2 \leq 1 - \frac{n-1}{N-1} = \frac{N-n}{N-1}.
$$
So $x_2 \geq \frac{1}{N-1},$ and so $x_j \geq \frac{1}{N-1}$ for $j=2,\ldots,N-1.$
It follows that $x_1 \leq 1 - \frac{N-1}{N-1} = 0,$ which gives a contradiction.
So there are no vertices with $1,2\notin J({\bf x}).$

{\bf Case 3.} $1 \notin J({\bf x})$ and $k \notin J({\bf x})$ for some $2<k\leq n+1.$
Again, since $\vert J({\bf x})\vert \geq N-1$ we must have $j \in J({\bf x})$ for $j \neq 1,k,$ so
\begin{equation}
\label{eq31}
0 < x_1 = x_2 = \cdots = x_{k-1} < x_k = x_{k+1} = \cdots = x_{n} = x_{n+1} = \cdots = x_{N-1} = x_N,
\end{equation}
\begin{equation}
\label{eq32}
x_1+x_2+\cdots+x_n=\frac{n-1}{N-1},
\end{equation}
and
\begin{equation}
\label{eq33}
x_1+x_2+\cdots+x_N=1.
\end{equation}

Using (\ref{eq31}) and (\ref{eq32}) we see that
$$
(k-1)x_1 + (n-(k-1)) x_k = \frac{n-1}{N-1},
$$
and from (\ref{eq31}) and (\ref{eq33}) we have
$$
(k-1)x_1 + (N-(k-1)) x_k = 1.
$$
Solving these two equations leads to $x_k = \frac{1}{N-1}$ and $x_1 = \frac{k-2}{(k-1)(N-1)}.$ Since $x_1
\leq x_k$ we conclude that
taking
$$
x_k = x_{k+1} = \cdots = x_N = \frac{1}{N-1},
$$
and
$$
x_1 = x_2 = \cdots = x_{k-1} = \frac{k-2}{(k-1)(N-1)}.
$$
leads to a point satisfying all of the inequalities. It is then straightforward to check the the rank
condition
so this is a vertex,
and this point corresponds to ${\bf w}^{(k-1)}.$

{\bf Case 4.} $1 \notin J({\bf x})$ and $k \notin J({\bf x})$ for some $n+1<k\leq N.$
We must have $j \in J({\bf x})$ for $j \neq 1,k,$ so
\begin{equation}
0 < x_1 = x_2 = \cdots = x_n = \cdots = x_{k-1} < x_k = x_{k+1} = \cdots = x_{N-1} = x_N,
\end{equation}
and in addition, (\ref{eq32}), (\ref{eq33}) hold.
Consequently
$$
nx_1 = \frac{n-1}{N-1},
$$
and
$$
(k-1)x_1 + (N-(k-1)) x_k = 1.
$$
Solving these two equations leads to $x_1 = \frac{n-1}{n(N-1)}$ and $x_k =
\frac{Nn-kn+k-1}{n(N-1)(N-k+1)}.$
It is straightforward to check that $x_1 \leq x_k,$ so all we conclude that, taking
$$
x_k = x_{k+1} = \cdots = x_N = \frac{Nn-kn+k-1}{n(N-1)(N-k+1)},
$$
and
$$
x_1 = x_2 = \cdots = x_{k-1} = \frac{n-1}{n(N-1)},
$$
leads to a point satisfying all of the inequalities.
It is then straightforward to check the the rank condition so gives a vertex point, and it is of the form
${\bf w}^{(k-1)}.$

{\bf Case 5.} $1 \notin J({\bf x})$ and $j \in J({\bf x})$ for $1\leq j \leq N.$
In this case, $0<x_1=x_j,$ for $j=2,\ldots,N.$ Thus, $x_j = \frac{1}{N}$ for $j=1,\ldots,N.$
Since $n<N$ it follows easily that inequality $I_{N+1}$ is strict, so $N+1 \notin J$
The rank condition is easily checked, so we obtain the vertex ${\bf w}^{(N)}.$
$~\Box$

{\bf Proof of Theorem 4.}
Given ${\bf x} = (x_1,\ldots,x_N)$ and a permutation ${\bf \pi} = (\pi_1,\ldots,\pi_N)$ of $(1,\ldots,N)$ we define
${\bf \pi}({\bf x}) = (x_{\pi_1},\ldots,x_{\pi_N}).$ We will refer to this as a permutation of ${\bf x}.$
Let ${\cal E}_{\bf \pi}$ denote the set of vertices of ${\cal O}_{\bf \pi}\cap {\cal T}_{N,n},$
and let ${\cal E} = \bigcup_{{\bf \pi}} {\cal E}_{\bf \pi}.$
Using Theorem 7, ${\cal E}_{\bf\pi}$ consists of those points of the form
${\bf\pi}({\bf w}^{(j)}),$ for $j=1,\ldots,N.$
Since ${\cal T}_{N,n}$ is the union of polytopes ${\cal O}_\pi \cap {\cal T}_{N,n}$ every vertex of
${\cal T}_{N,n}$ lies in ${\cal E}.$
We proceed to identify certain points in ${\cal E}$ as averages of permutations
of the points ${\bf w}^{(1)}$ and ${\bf w}^{(N-1)},$ which shows they cannot be vertices of ${\cal T}_{N,n}.$

Define points

$$
\begin{array}{c}
{\bf u}^{(1)} = (\underbrace{a_1,b_1,\ldots,b_1}_{n},\underbrace{b_1,\ldots,b_1}_{N-n})\\
{\bf u}^{(2)} = (\underbrace{b_1,a_1,b_1,\ldots,b_1}_{n},\underbrace{b_1,\ldots,b_1}_{N-n})\\
{\bf u}^{(3)} = (\underbrace{b_1,b_1,a_1,b_1,\ldots,b_1}_{n},\underbrace{b_1,\ldots,b_1}_{N-n})\\
\vdots \\
{\bf u}^{(n-1)} = (\underbrace{b_1,\ldots,b_1,a_1,b_1}_{n},\underbrace{b_1,\ldots,b_1}_{N-n})\\
{\bf u}^{(n)} = (\underbrace{b_1,b_1,\ldots,b_1,a_1}_{n},\underbrace{b_1,\ldots,b_1}_{N-n})\\
\end{array}
$$

each of which is a permutation of ${\bf w}^{(1)}.$
It follows that
$\sum_{i=1}^n u^{(k)}_i = \frac{n-1}{N-1},\textrm{ for } k=1,\ldots,n.$

For fixed  $j$ with $1 \leq j\leq n$ define
$$
{\bf u} = \frac{1}{j} \sum_{k=1}^j {\bf u}^{(k)}.
$$
Then we have
$$
{\bf u} = (\underbrace{a,\cdots,a}_{j},\underbrace{b,\cdots,b}_{N-j})
$$
for some values of $a$ and $b.$ In addition, since, for $k=1,\ldots,j$ we have
$$
\sum_{i=1}^n u^{(k)}_i = \frac{n-1}{N-1},
$$
and
$$
\sum_{i=1}^N u^{(k)}_i = 1,
$$
it follows that
$$
\sum_{i=1}^n u_i = \frac{n-1}{N-1},
$$
and
$$
\sum_{i=1}^N u_i = 1.
$$
It follows (see Remark 1) that $a = a_j$ and $b = b_j,$ so ${\bf u} = {\bf w}^{(j)}.$
Thus, ${\bf w}^{(j)}$ is not a vertex of ${\cal T}_{N,n}.$
By symmetry, the same is true of ${\bf \pi}({\bf w}^{(j)})$ for all permutations $\pi.$

In addition, we have ${\bf w}^{(N)} = (1/N,\ldots,1/N)$ and this is the average of the 
${\bf u}^{(k)},~k=1,\ldots,N,$
so ${\bf w}^{(N)}$ is not a vertex of ${\cal T}_{N,n}.$ Nor of course, are its permutations.

Next, define points
$$
\begin{array}{c}
{\bf v}^{(n+1)} = (\underbrace{a_{N-1},\cdots,a_{N-1}}_{n},\underbrace{b_{N-1},a_{N-1}\ldots,a_{N-1}}_{N-n})\\
{\bf v}^{(n+2)} =
(\underbrace{a_{N-1},\cdots,a_{N-1}}_{n},\underbrace{a_{N-1},b_{N-1},a_{N-1}\ldots,a_{N-1}}_{N-n})\\
{\bf v}^{(n+3)} =
(\underbrace{a_{N-1},\cdots,a_{N-1}}_{n},\underbrace{a_{N-1},a_{N-1},b_{N-1},a_{N-1}\ldots,a_{N-1}}_{N-n})\\
\vdots \\
{\bf v}^{(N-1)} =
(\underbrace{a_{N-1},\ldots,a_{N-1}}_{n},\underbrace{a_{N-1},\ldots,a_{N-1},b_{N-1},a_{N-1}}_{N-n})\\
{\bf v}^{(N)} = (\underbrace{a_{N-1},\ldots,a_{N-1}}_{n},\underbrace{a_{N-1},\ldots,a_{N-1},b_{N-1}}_{N-n})\\
\end{array},
$$
each of which is a permutation of ${\bf w}^{(N-1)}.$
Since $na_{N-1} = \frac{n-1}{N-1}$ we see that
$\sum_{i=1}^n v^{(k)}_i = \frac{n-1}{N-1},\textrm{ for } k=n+1,\ldots,N.$

For fixed  $j$ with $n \leq j\leq N-1$ define
$$
{\bf v} = \frac{1}{N-j} \sum_{k=j+1}^N {\bf u}^{(k)}.
$$
Then 
$$
{\bf v} = (\underbrace{a,\cdots,a}_{j},\underbrace{b,\cdots,b}_{N-j})
$$
for some values of $a$ and $b.$ In addition, since, for $k=j+1,\ldots,N$ we have
$$
\sum_{i=1}^n v^{(k)}_i = \frac{n-1}{N-1},
$$
and
$$
\sum_{i=1}^N v^{(k)}_i = 1,
$$
it follows that
$$
\sum_{i=1}^n v_i = \frac{n-1}{N-1},
$$
and
$$
\sum_{i=1}^N v_i = 1.
$$
It follows that $a = a_j$ and $b = b_j,$ so ${\bf u} = {\bf w}^{(j)}.$
Thus, ${\bf w}^{(j)}$ is not a vertex of ${\cal T}_{N,n}.$
By symmetry, the same is true of ${\bf\pi}({\bf w})^{(j)}$ for all permutations ${\bf\pi}.$

This proves that the permutations of ${\bf w}^{(1)},$ i.e. points of the form ${\bf p}^{(i,0)},$ and
the permutations of ${\bf w}^{(N-1)},$ i.e. points of the form ${\bf p}^{(i,\frac{1}{n})},$ for $i=1,\ldots,N,$
are the only possibilities for vertices of ${\cal T}_{N,n}.$
Note that
$p^{(k,0)}_k = 0,$ but
$$
p^{(j,0)}_k = \frac{1}{N-1} > 0, \textrm{ for } j\neq k,
$$
and
$$
p^{(j,1/n)}_k = \frac{n-1}{n(N-1)} > 0, \textrm{ for } j=1,\ldots,N.
$$
Since any convex combination of the points ${\bf p}^{(j,0)}, ~j\neq k,$ and ${\bf p}^{(j,1/n)},~j=1,\ldots,N,$
must have a positive $k^{\textrm{th}}$ coordinate, no convex combination can equal ${\bf p}^{(k,0)},$
and we conclude that ${\bf p}^{(k,0)}$ is a vertex of ${\cal T}_{N,n}$ for $k=1,\ldots,N.$

If $n<N-1$ then
$p^{(k,1/n)}_k = \frac{1}{n},$
but
$$
p^{(j,1/n)}_k = \frac{n-1}{n(N-1)} < \frac{1}{n}, \textrm{ for } j\neq k,
$$
and in addition,
$$
p^{(j,0)}_k = \frac{1}{N-1} < \frac{1}{n}, \textrm{ for } j=1,\ldots,N.
$$
Thus, any convex combination of ${\bf p}^{(j,1/n)}, ~j\neq k,$ and ${\bf p}^{(j,0)},~j=1,\ldots,N$
must have as its $k^{\textrm{th}}$ coordinate, a value less than $\frac{1}{n},$
so ${\bf p}^{(k,1/n)}$ cannot be expressed as such a convex combination, and we conclude that
each ${\bf p}^{(k,1/n)}$ is a vertex of ${\cal T}_{N,n}$ for $k=1,\ldots,N.$

Finally, in case $n=N-1,$ it is straightforward to check that
$$
{\bf p}^{(k,1/n)} = \frac{1}{N-1} \sum_{j=1,j\neq k}^N {\bf p}^{(j,0)}
$$
so none of the ${\bf p}^{(k,1/n)}$ is a vertex.$~~\Box$

{\bf Proof of Theorem 5.} ${\cal T}_{N,n}$ is defined to be those points ${\bf x} = (x_1,\ldots,x_N)$
satisfying
(\ref{faceconstraints}) together with
\begin{equation}
\label{equalityconstraint}
\sum_{i=1}^N x_i = 1,
\end{equation}
\begin{equation}
\label{nonnegativityconstraint}
x_i \geq 0, \textrm{ for } i=1,\ldots,N.
\end{equation}
We proceed to demonstrate that (\ref{nonnegativityconstraint}) is a consequence of (\ref{faceconstraints})
and (\ref{equalityconstraint}).

Assume (\ref{faceconstraints}) and (\ref{equalityconstraint}) and fix $i \in \{ 1,\ldots,N\}.$  Then
summing
(\ref{faceconstraints}) over subsets $F \in {\cal P}_N(n)$ containing $i$ we see that
\begin{equation}
\label{sumoffaceconstraints}
\sum_{F\in {\cal P}_N(n)~:~ i\in F} \sum_{j\in F} x_j \geq {N-1 \choose n-1}\frac{n-1}{N-1} = {N-2 \choose
n-2}.
\end{equation}
On the other hand, consider the sum on the left-hand side of (\ref{sumoffaceconstraints}).
In the sum, $x_i$ appears in every term, hence ${N-1 \choose n-1}$ times, while for $j \neq i$ the number
of times $x_j$ appears is
the number of subsets of $\{1,\ldots,N\}$ of size $n-1$ that contain $j$ and do not contain $i,$ which is
${N-2 \choose n-2}.$
Thus, the left-hand side takes the form
\begin{eqnarray*}
{N-1 \choose n-1}x_i + {N-2 \choose n-2}\sum_{1\leq j \leq N, j\neq i} x_j  &=&
\left({N-1 \choose n-1}-{N-2 \choose n-2}\right)x_i + {N-2 \choose n-2}\sum_{j=1}^n x_j \\
&=& {N-2 \choose n-1} x_i + {N-2 \choose n-2}.
\end{eqnarray*}
So we can conclude that
$$
{N-2 \choose n-1} x_i + {N-2 \choose n-2} \geq {N-2 \choose n-2},
$$
so $x_i \geq 0.$

If any of the inequality constraints in (\ref{faceconstraints}) fails to define a facet of ${\cal
T}_{N,n}$ then, by symmetry,
none of them would define a facet, which is clearly not possible, so each must define a facet.

Next, consider the question of which vertices lie in which facets. The
vertices ${\bf x}$ in the facet defined by (\ref{faceconstraints}) for a given $F \in {\cal P}_N(n),$ are those
for which
we have
$$
\sum_{i\in F} x_i = \frac{n-1}{N-1}.
$$
Clearly, this equality holds for ${\bf x} = {\bf p}^{(i,0)}$ if $i\in F,$ and fails for $i\notin F.$
In addition, for $i\notin F$ we have
$$
\sum_{j\in F} p^{(i,1/n)}_j = n\frac{n-1}{n(N-1)} = \frac{n-1}{N-1},
$$
while for $i \in F,$ it is the case that
$$
\sum_{j\in F} p^{(i,1/n)}_j = \frac{1}{n} + (n-1)\frac{n-1}{n(N-1)}
$$
$$
= \frac{n-1}{N-1}
\left[ \frac{1}{n} \left\{
\frac{N-1}{n-1}+n-1\right\}\right] >
\frac{n-1}{N-1}
\left[ \frac{1}{n} \left\{
1+n-1\right\}\right] = \frac{n-1}{N-1}.
$$
$~~\Box$

{\bf Proof of Lemma 2.}
We make use of some easily proven facts.
each of which can proven by taking the log, replacing the log by an integral, and using an elementary inequality.
First, for $1 \leq r \leq s$ we have
$$
e^{-r - \frac{r^2}{s-r}} \leq (1-r/s)^s \leq e^{-r}.
$$
Second, for $x\in (-1,\infty)$
$$
\frac{x}{1+x} \leq \log(1+x) \leq x.
$$

Applying the improved Stirling approximation \cite{robbins} we have
\begin{equation}
\label{Stirling1}
\frac{s!}{(s-r)!} \leq
\left(\frac{s}{s-r}\right)^s (s-r)^r \left(\frac{s}{s-r}\right)^\frac{1}{2}
e^{\frac{1}{12s}-\frac{1}{12(s-r)+1}}e^{-r}.
\end{equation}
We proceed to upper bound each of the first four terms in (\ref{Stirling1}).
First, we have
\begin{equation}
\label{lemmaineq1}
\left( \frac{s}{s-r} \right)^s = \left( 1- r/s \right)^{-s} = e^{-s \log(1-r/s)}\leq e^{-s \frac{-r/s}{1-r/s}} =
e^{r+\frac{r^2}{s-r}}.
\end{equation}
For the second term, we have
\begin{equation}
\label{lemmaineq2}
(s-r)^r = s^r (1-r/s)^r = s^r e^{r \log(1-r/s)}\leq s^r e^{-\frac{r^2}{s}}.
\end{equation}
For the third term
\begin{equation}
\label{lemmaineq3}
\left( \frac{s}{s-r}\right)^\frac{1}{2} = e^{-\frac{1}{2} \log(1-r/s)} \leq e^{-\frac{1}{2} \frac{-r/s}{1-r/s}}
=  e^{\frac{r}{2(s-r)}}.
\end{equation}
Finally,
\begin{equation}
\label{lemmaineq4}
e^{\frac{1}{12s}-\frac{1}{12(s-r)+1}} \leq e^{\frac{1}{12s}-\frac{1}{12s+1}} = e^{\frac{1}{144s^2}}
\end{equation}
Combining (\ref{lemmaineq1}), (\ref{lemmaineq2}), (\ref{lemmaineq3}), and (\ref{lemmaineq4})
gives the claimed lower bound for $\frac{(s-r)!}{s!}.$

For the upper bound, again, using the improved Stirling approximation, we have
\begin{equation}
\label{Stirling2}
\frac{(s-r)!}{s!} \leq
\left(\frac{s-r}{s}\right)^s (s-r)^{-r} \left(\frac{s-r}{s}\right)^\frac{1}{2}
e^{\frac{1}{12(s-r)}-\frac{1}{12s+1}}e^{r}.
\end{equation}
Again, we proceed to upper bound each of the first four terms in (\ref{Stirling2}).
For the first term
\begin{equation}
\label{lemmaineq5}
\left( \frac{s-r}{s} \right)^s = \left( 1- r/s \right)^{s} = e^{s \log(1-r/s)}\leq e^{-r}.
\end{equation}
For the second term, we have
\begin{equation}
\label{lemmaineq6}
(s-r)^{-r} = s^{-r} (1-r/s)^{-r} = s^{-r} e^{-r \log(1-r/s)}\leq s^{-r} e^{-r\frac{-r/s}{1-r/s}} = s^{-r} e^{\frac{r^2}{s-r}}
\end{equation}
For the third term
\begin{equation}
\label{lemmaineq7}
\left( \frac{s-r}{s}\right)^\frac{1}{2} = e^{\frac{1}{2} \log(1-r/s)} \leq e^{-\frac{r}{2s}}.
\end{equation}
Finally,
\begin{equation}
\label{lemmaineq8}
e^{\frac{1}{12(s-r)}-\frac{1}{12s+1}} \leq e^{\frac{1}{12(s-r)}}.
\end{equation}
Combining (\ref{lemmaineq5}), (\ref{lemmaineq6}), (\ref{lemmaineq7}), and (\ref{lemmaineq8})
gives the claimed upper bound for $\frac{(s-r)!}{s!}.$
$~~\Box$

\bibliography{mybib}
\bibliographystyle{plain}

\end{document}